\documentclass[twoside,a4paper,reqno,12pt]{amsart}
\usepackage{lmodern}
\usepackage{amsfonts, amsbsy, amsmath, amssymb, latexsym}
\usepackage{mathrsfs,array}
\usepackage[top=1 in,right=1 in,bottom=1 in,left=1 in]{geometry}
\usepackage[utf8]{inputenc}
\usepackage{xcolor}
\usepackage[pagebackref,colorlinks,citecolor=blue,linkcolor=magenta]{hyperref}
\usepackage{tikz}
\usepackage[english]{babel}
\usepackage{etoolbox}
\usepackage{faktor}
\usepackage{graphicx}
\usepackage{mathtools}
\usepackage{floatflt}
\usepackage{multicol}
\usepackage{lipsum}
\usepackage{diagbox}
\usepackage{wrapfig}
 \usepackage{relsize}
 \usepackage{booktabs}
\usepackage{pgfplots}
\usepackage{pgfplotstable}
\usepackage{enumerate}
\usepackage{blindtext}
\usepackage{amssymb}
\usepackage{amsthm}
\usepackage{amsmath}
\usepackage{mathabx}
\usepackage{siunitx}
\usepackage{pdflscape}
\newcommand\dyckpath[5]{
  \begin{scope}[local bounding box=#4]
    \fill[white]  (#1) rectangle +(#2,#2);
    \fill[blue!25!white] (#1) foreach \dir in {#3}{-- ++(\dir*90:1)} |- (#1);
    
    \path[fill] (#1) foreach \num [count=\i from 0] in {#5}{ +(\i,0) node[anchor=north]{\num} \ifnum\i>#2 circle (1pt) \fi};
    
    \draw[help lines] (#1) grid +(#2,#2);
    
    \draw[line width=2pt] (#1) foreach \dir in {#3}{ -- ++(\dir*90:1)};
    
    \draw[dashed, blue] (#1) -- +(#2,#2);
    
  \end{scope}
}
\DeclareMathAccent{\widecheck}{0}{mathx}{"71}

\newcommand{\id}{\operatorname{id}}

\newcommand\Des{\operatorname{Des}} 
\newcommand\des{\operatorname{des}} 
\newcommand\cC{\mathcal{C}}
\newcommand\cL{\mathcal{L}}

\usetikzlibrary{lindenmayersystems,arrows.meta}

\newcommand\commentout[1]{}

\usepackage{array,tabularx}

\usepackage[english]{babel}

\theoremstyle{plain}
\newtheorem{theorem}{Theorem}[section]
\newtheorem{proposition}[theorem]{Proposition}

\newtheorem{corollary}[theorem]{Corollary}
\newtheorem{conj}[theorem]{Conjecture}

\theoremstyle{definition}
\newtheorem{definition}[theorem]{Definition}
\newtheorem{example}[theorem]{Example}
\newtheorem{question}[theorem]{Question}

\theoremstyle{remark}
\newtheorem{remk}[theorem]{Remark}

\newtheorem*{conj*}{Conjecture}

\usepackage{graphicx} 
\usepackage{amsmath}

\title{Canon Permutation Posets}

\author{Matthias Beck}
\address{Department of Mathematics\\
         San Francisco State University\\
         San Francisco, CA 94132\\
         U.S.A.}
\email{mattbeck@sfsu.edu}

\author{Danai Deligeorgaki}
\address{Department of Mathematics, Universitat de Barcelona, Barcelona, Spain}
\email{deligeorgaki@ub.edu}

\date{13 September 2026}



\begin{document}

\maketitle

\begin{abstract}
A permutation of the multiset $\{1^m,2^m,\dots,n^m\}$ is a {\em canon permutation} if
the subsequence formed by the $j$th copy of each element of $[n]:=\{1,2,\dots,n\}$ is identical for
all $j\in[m]$. Canon permutations were introduced by Elizalde and are motivated by pattern-avoiding 
concepts such as (quasi-)Stirling permutations. He proved that the descent polynomial of canon
permutations exhibits a surprising product structure; as a further consequence, it is palindromic.
Our goal is to understand canon permutations from the viewpoint of Stanley's $(P,\omega)$-partitions,
along the way generalizing Elizalde's definition and results. We start with a labeled poset $P$ and
extend it in a natural way to \emph{canon labelings} of the product poset $P \times [n]$. The resulting descent
polynomial has a product structure which arises naturally from the theory of $(P,\omega)$-partitions and simplifies existing proofs. When
$P$ is graded, this theory also implies palindromicity. We include results on weak descent
polynomials, an amphibian construction between canon permutations and multiset permutations, giving rise to
\emph{dissonant canon permutations}, as well as
$\gamma$-positivity and interpretations of descent polynomials of canon permutations.
\end{abstract}


\section{Introduction}

A permutation $\pi$ of the multiset $\{1^m,2^m,\dots,n^m\}$ is a {\em canon permutation} if
the subsequence formed by the $j$th copy of each element of $[n]:=\{1,2,
\dots,n\}$ is identical for all $j\in[m]$. For instance, $223143213144$ is a canon
permutation of $\{1^3,2^3,3^3,4^3\}$, with the sequence formed by the $j$th copy being $2314$.
Canon permutations were introduced by Elizalde~\cite{sergi} and are motivated by
pattern-avoiding permutations of the multiset $\{1,1,2,2,\dots,n,n\}$, such as Stirling
\cite{GESSEL197824} and quasi-Stirling \cite{Archer2018PatternRQ} permutations. 
For $m = 2$, canon permutations are sometimes called \emph{nonnesting permutations}; they
are precisely the permutations of the multiset $\{1,1,2,2,\dots,n,n\}$ avoiding the patterns
$1221$ and $2112$.
Following the notation
in~\cite{sergi}, we denote by $\cC_n^{m,\sigma}$ the set of canon permutations whose $j$th copy of each entry forms the permutation $\sigma\in S_{n}$ and by $\cC_n^{m}$ the set of all canon permutations for given $m,n$. Here, $S_{n}$ denotes the set of permutations of $[n]$.

As usual, we call $j$ a \emph{descent} of $\pi$ if $\pi(j+1) < \pi(j)$ and denote by
$\des(\pi)$ the number of descents of $\pi$. We further denote by $\operatorname{wdes}(\pi)$ the
number of \emph{weak descents} of $\pi$, i.e., the number of positions $j$ such that
$\pi(j+1) \leq  \pi(j)$. Elizalde's work 
centered on
understanding the distribution of descents of canon permutations by studying the
\emph{descent polynomial}
\[C^{m}_n(x):=\sum_{\pi \in \cC_n^{m}}x^{\des(\pi)} . \]

Let $A_n(x)$ be an Eulerian polynomial (i.e., the descent polynomial of
permutations on $[n]$) and $N_n(x)$ a Narayana polynomial, enumerating high peaks in
Dyck paths (we give a precise definition later).
Both of these polynomials are \emph{palindromic}, that is, their coefficient
sequences are symmetric.
Elizalde found that the descent polynomial $C^2_n(x)$ has the following
surprisingly simple structure~\cite[Theorem~2.1]{sergi}.

\begin{theorem}[Elizalde \cite{sergi}]\label{thm: sergi}
For $n \geq 1$, we have $C^2_n(x) = A_n(x) \, N_n(x)$.
In particular, $C^2_n(x)$ is palindromic.
\end{theorem}

As Elizalde noted, the palindromicity of $C^2_n(x)$ is a priori unexpected.
He proved a bivariate generalization of Theorem~\ref{thm: sergi}, involving also the number of the so-called \textit{plateaux} of permutations. 
Another consequence of the bivariate extension of Theorem~\ref{thm: sergi} is that the number of
canon permutations with $r$ weak descents is the same as that for $2n-r$ descents. (We will see an explanation of this palindromicity in Proposition~\ref{prop:weakdes}.)
In a follow-up paper~\cite{sergi2}, Elizalde extended Theorem~\ref{thm: sergi} from $m = 2$
to the general case. Here the Narayana polynomial is generalized so that the role of Dyck
paths is replaced by standard Young tableaux of rectangular shape. Elizalde's proofs are
bijective, and the resulting descent polynomial $C^m_n(x)$ is, again surprisingly,
palindromic, due to the (older) fact that the generalized Narayana polynomials are
palindromic~\cite{sulanke}. 

Our goal is to understand Theorem~\ref{thm: sergi} and its extensions from the viewpoint
of Stanley's $(P,\omega)$-partitions~\cite{stanleythesis,ec1}, along the way generalizing
Elizalde's definition and results. Elizalde asked for a more direct explanation of the structure
of his descent polynomials, and this viewpoint delivers such an explanation, simplifying the
existing proofs along the way. We describe our \emph{ansatz} next.
Consider a poset $P$ with $m$ elements 
and a labeling $\omega$ (i.e., a bijection $\omega: P \to [m]$). We extend this labeling, via a
given $\sigma \in S_n$, to the {\em canon labeling} $\omega \times \sigma$ of the poset $P \times [n]$
(here we think of $[n]$ as a poset---a chain) defined by
\[
  (\omega \times \sigma)(p, j) := \omega(p) + (\sigma(j) - 1)m \, .
\]
A \emph{linear extension} of $P \times [n]$ is defined, as usual, as an
order-preserving bijection $\pi : [mn] \to P \times [n]$; writing $\pi$ as a word in
terms of the canon labeling $\omega \times \sigma$ of $P \times [n]$, it is a short step to view $\pi$ as a
permutation of the multiset $\{1^m,2^m,\dots,n^m\}$. 
Let $\cC_n^{ P, \omega \times \sigma }$ consist of all such multiset permutations, and define
\[
  \cC_n^{ P, \omega } := \left\{ \cC_n^{ P, \omega \times \sigma } : \, \sigma \in S_n \right\} 
\]
which we call the set of all \emph{canon permutations} of the labeled poset $(P,\omega)$.
The classical canon permutations  $\cC_n^m$ stem from the case in which $P=[m]$ is a chain and $\omega$ is
the identity (i.e., $[m]$ is \emph{naturally labeled}). 
For instance, starting with the $2$-element chain poset $P=[2]$ and the identity permutation $\id\in S_4$,
the product poset $[2]\times [4]$ is  shown in Figure~\ref{fig:aaaaa}; the canon labeling $\omega \times \sigma$ of $[2]\times [4]$ is identified with a natural labeling with entries in $[8]$. 
Here, $\cC_4^{ P, \omega }$ can be realized as the set $\cC_4^2$ of canon permutations of $\{1,1,2,2,3,3,4,4\}$.

\begin{figure}[ht]
    \centering
\begin{tikzpicture}[main/.style = {}] 
{\color{magenta}
\node[main] (1) {$1$}; 
\node[main] (2) [above of=1] {$2$};
\node[main] (3) [right of=2] {$3$};
\node[main] (4) [above of=3] {$4$};
\node[main] (5) [right of=4] {$5$};
\node[main] (6) [above of=5] {$6$};
\node[main] (7) [right of=6] {$7$};
\node[main] (8) [above of=7] {$8$};
\draw (2) -- (1);
\draw (3) -- (4);
\draw (5) -- (6);
\draw (7) -- (8);
\draw (3) -- (1);
\draw (2) -- (4);
\draw (5) -- (7);
\draw (6) -- (8);
\draw (3) -- (5);
\draw (6) -- (4);}
\end{tikzpicture} 
\hspace{2 cm}
\begin{tikzpicture}[main/.style = {}] 
{\color{blue}
\node[main] (4) {$1$}; 
\node[main] (2) [above of=1] {$2$};
\node[main] (3) [right of=2] {$5$};
\node[main] (4) [above of=3] {$6$};
\node[main] (5) [right of=4] {$3$};
\node[main] (6) [above of=5] {$4$};
\node[main] (7) [right of=6] {$7$};
\node[main] (8) [above of=7] {$8$};
\draw (2) -- (1);
\draw (3) -- (4);
\draw (5) -- (6);
\draw (7) -- (8);
\draw (3) -- (1);
\draw (2) -- (4);
\draw (5) -- (7);
\draw (6) -- (8);
\draw (3) -- (5);
\draw (6) -- (4);}
\end{tikzpicture} 
    \caption{The poset $[m]\times [n]$ for $m = 2$, $n = 4$ and $\sigma = \id \in S_4$ \textcolor{magenta}{(left)}, or  $\sigma = 1324 \in S_4$ \textcolor{blue}{(right)}. On the \textcolor{magenta}{left}, the linear extension $1 3 2 4 5 7 6 8$ is identified with
the canon permutation $1  2 1 2 3 4 3 4$. On the \textcolor{blue}{right}, $1 5 2 6 3 7 4 8$
corresponds to $1 3 1 3 2 4 2 4$.}
    \label{fig:aaaaa}
\end{figure}

Finally, we define the \emph{canon polynomial} 
\[C^{P,\omega}_n(x):=\sum_{\pi \in \mathcal{C}_n^{P,\omega}}x^{\des(\pi)} . \]
Again, we note the special case $C_n^m(x) = C_n^{ [m], \id } (x)$.
Indeed, in all but one of our applications, we will use a natural labeling for~$\omega$.
The one exception (at least in this current work) is captured by
the \emph{weak descent polynomial} $ C_n^{ P, {w_0} } (x)$, where the labeling ${w_0}$ is reverse
natural. In the case $P=[m]$, this corresponds to the polynomial
enumerating weak descents in canon permutations of $\{1^m,2^m, \dots, n^m\}$. We will see
that it is a translate of $C_n^{P,\omega}(x)$ for the posets that we study
(Proposition~\ref{prop:weakdes} below). 
Our generalization of Theorem~\ref{thm: sergi} and its extensions is as follows.

\begin{theorem}\label{thm:main}
Let $P$ be a poset with a natural labeling $\omega$. Then $C_n^{P,\omega}(x) = A_n(x) \, h^*_{P \times [n]}(x)$.
Furthermore, if $P$ is graded then $C_n^{P,\omega}(x)$ is palindromic.
\end{theorem}

Here $P$ being \emph{graded} means that all maximal chains have the same length, 
and $h^*_{P \times [n]}(x)$ is the numerator of the rational generating function of the
order polynomial of $P \times [n]$, which is in a sense the descent polynomial of $P
\times [n]$; we will give detailed definitions in Section~\ref{sec:ordpollabelposets}.
In fact, Theorem \ref{thm:main} is a special case of a more general result (Theorem~\ref{thm:
multiplication} below).

In Section~\ref{sec:narayana}, we give a poset model whose descent polynomial is the Narayana
polynomial. The underlying posets turn out to be structured such that all labeled versions of them, in the
sense of $(P,\omega)$-partitions, have similar descent polynomials, from which Theorem~\ref{thm:
sergi} follows in a few short steps, as we outline in Section~\ref{sec:ordpollabelposets}.
The philosophy of our \emph{ansatz} is that we consider the poset from
Section~\ref{sec:narayana}, whose descent polynomial is the Narayana polynomial,
and then sum $(P,\omega)$-descent polynomials over certain labelings $\omega$ of this fixed
poset, which gives rise to the product structure exhibited in Theorem~\ref{thm: sergi}.
Our proof generalizes immediately to the canon polynomial $C^m_n(x)$ for general~$m$.
The palindromicity of canon polynomials follows organically from the structure of the
involved posets;
indeed, we give three direct explanations for this palindromicity (Theorem~\ref{thm: multiplication}, Proposition~\ref{prop:oneposet}, and Theorem~\ref{them: subposet} below).

In Section~\ref{sec:amphibians}, we extend our results to 
subposets of $P\times [n]$ with some relations missing between different
copies of $P$, i.e., those of the form $(p,j) \prec (p,j+1)$. 
These subposets give intermediary descent polynomials situated between those for canon permutations
and multiset permutations which we call \textit{dissonant canon permutations}. We show that palindromicity extends to this intermediary class.

A distributional property that is stronger than palindromicity (and unimodality \cite{pettersurvey}) is $\gamma$-positivity: a palindromic polynomial of degree $d$ is {\em $\gamma$-positive} if its coefficients are non-negative when expressed in the \emph{$\gamma$-basis}
$ \{x^i(x + 1)^{d - 2i} : 0 \leq i \leq \lfloor d/2\rfloor\}$ (see for instance \cite{athanasiadis2017gamma}).
In Section~\ref{sec:amphibians}, we further  conjecture that the intermediate class of polynomials
have $\gamma$-positive descent polynomials and ask
for a $\gamma$-coefficient interpretation, which would recover $\gamma$-positivity exhibited by both canon permutations and multiset permutations. In Section~\ref{sec: canon positivity}, we show that our viewpoint implies that the descent 
polynomials of canon permutations are $\gamma$-positive and give a combinatorial interpretation for their coefficients,
based on a result of Br\"and\'en.


\section{Narayana Polynomials as Descent Polynomials}\label{sec:narayana}

We start, as a warm-up of sorts, by realizing the Narayana polynomials as descent
polynomials of the posets $[2] \times [n]$.
As usual, the set $D_n$ of \emph{Dyck paths} consists of all lattice paths from $(0, 0)$ to
$(n, n)$ with steps $e := (1, 0)$ and $n := (0, 1)$ that do not go above the diagonal $y = x$. 
A \emph{peak} in a Dyck path is an occurrence of two adjacent
steps $en$. 
A peak is called a \emph{high peak} if these steps do not touch the diagonal.
We denote the number of high peaks of $D \in D_n$ by $\text{hpea}(D)$. For example, the Dyck path from $(0,0)$ to $(4,4)$ appearing in Figure \ref{fig: labeled poset} has three peaks out of which one is a high peak. We refer to the polynomial 
 $N_n(x):=\sum_{D \in D_n} x^{\text{\rm hpea}(D)}$ as a {\em Narayana polynomial}. 
(Our convention is slightly nontraditional, as the \emph{Narayana numbers} are the coefficients of $x
\, N_n(x)$.)

\begin{figure}[ht]
    \centering
\begin{tikzpicture}
{\color{blue}[main/.style = {}] 
\node (1) {$1$}; 
\node (2) [above of=1] {$2$};
\node (3) [right of=2] {$3$};
\node (4) [above of=3] {$4$};
\node (5) [right of=4] {$5$};
\node (6) [above of=5] {$6$};
\node (7) [right of=6] {$7$};
\node (8) [above of=7] {$8$};
\draw (2) -- (1);
\draw (3) -- (4);
\draw (5) -- (6);
\draw (7) -- (8);
\draw (3) -- (1);
\draw (2) -- (4);
\draw (5) -- (7);
\draw (6) -- (8);
\draw (3) -- (5);
\draw (6) -- (4);
}
\end{tikzpicture} 
\hspace{2.5 cm}
{\color{blue}
\begin{tikzpicture}[scale=0.9]

  \dyckpath{0,0}{4}{0,1,0,0,1,1,0,1}{dyck4}{ };

\end{tikzpicture}}
    \caption{The poset $[2]
\times [4]$ naturally labeled \textcolor{blue}{(left)}. The Dyck path corresponding to the linear extension $12354678$ of $[2]
\times [4]$ \textcolor{blue}{(right)}.}
    \label{fig: labeled poset}
\end{figure}

 \begin{theorem}\label{thm: narayana bijection}
The linear extensions of the poset $[2]
\times [n]$ are in bijection with the Dyck paths in $D_n$. Furthermore, the descents of a
linear extension of $[2] \times [n]$, labeled via $\id \times \id$, are in bijection with the high peaks of the corresponding Dyck path. 
 \end{theorem}
 
 Considering the equivalence between linear extensions of $[2]\times [n]$ and canon
permutations of $\{1^2,2^2,...,n^2\}$, and therefore ${C}_n^{2}(x)={C}_n^{[2],\id}(x)$, one can also derive Theorem~\ref{thm: narayana bijection} from a bijection of Elizalde, constructed in \cite[Section~3.1]{sergi}.


\begin{proof}
We label $[2] \times [n]$ naturally via $\id \times \id$; see Figure~\ref{fig: labeled
poset} (left) as an example.
     We map a linear extension $\pi=\pi_1\cdots\pi_{2n}$ of $[2] \times [n]$ to the Dyck path $D=d_1d_2 \cdots d_{2n}$ where $d_i=e$, if $\pi_i$ is odd, and
$d_i=n$, otherwise. 
Since the order in which the odd (resp.\ even) numbers appear in a
linear extension of $([2] \times [n],\id \times \id)$ is fixed, the map is invertible,  giving
a bijection between the set of linear extensions of $[2] \times [n]$ and $D_n$.

     Now, let us discuss the correspondence between descents of a linear extension
$\pi=\pi_1\cdots\pi_{2n}$ of $([2] \times [n],\id \times \id)$ and high peaks of the corresponding Dyck
path $D=d_1\cdots d_{2n}$. 
We have that $i$ is a descent of $\pi$  ($2\leq i \leq 2n-2$ by
construction of $[2] \times [n]$) if and only if $\pi_i=2k+1$ for some $1\leq k \leq n-1$ and
$\pi_{i+1}=2k-2j$ for some $0\leq j \leq k-1$. 
This happens if and
only if the entries $\pi_1, \dots,\pi_{i-1}$ contain more odd numbers than even and
$\pi_i$ is odd but $\pi_{i+1}$ is even. That means exactly that the entries $d_1,
\dots, d_{i-1}$ of $D$ consist of more east than north steps, and that $d_i=e$ while $d_{i+1}=n$. 
By definition, that happens if and only if $d_id_{i+1}$ is a high peak. 
 \end{proof}
For example, the linear extension $12354678\in \mathcal{L}([2]\times [4], \id\times \id)$ corresponds to the Dyck path shown in Figure \ref{fig: labeled poset}.

\begin{corollary}\label{cor:h*asnarayana} 
Let $\cL([2] \times [n])$ denote the set of linear extensions of $[2] \times [n]$. Then
     \[\sum_{\pi \in \cL([2]\times [n])} x^{\des(\pi)}=N_n(x)\, . \]
 \end{corollary}

The analogy between descents of linear extensions and high peaks of Dyck paths extends to weak
descents and (ordinary) peaks, and the bijection discussed extends similarly.\footnote{We note that Corollary \ref{cor:h*asnarayana} also follows from work of Sulanke \cite{sulanke}.}

\begin{remk}
Viewing the descent polynomial of $[2] \times [n]$ as the Ehrhart $h^*$-polynomial of the order
polytope of $[2] \times [n]$, it is a curious fact that the convex hull of the origin and the positive simple roots
in $A_{ n-1 }$ has the same $h^*$-polynomial~\cite[Example~6]{braunsurvey}.
There is another connection between this poset and type-A roots (but it is unclear to us whether this is related):
the poset of order ideals of $[2] \times [n]$ naturally corresponds bijectively to the \emph{root poset} of $A_n$, given by the positive roots and the relation $\alpha \preceq \beta$ whenever $\beta - \alpha$ is a sum of simple positive roots.
\end{remk}


\section{Order Polynomials of Labeled Posets}\label{sec:ordpollabelposets} 

We now recall some fundamental definitions and results on labelled posets and their order polynomials
and generating functions~\cite{stanleythesis,ec1}.

\begin{definition}
Let $P$ be a poset of cardinality $m$ with a given labeling $\omega: P \to [m]$.
A \emph{$(P, \omega)$-partition} is a map $\sigma: P \rightarrow \mathbb{Z}_{ \ge 0 } $
satisfying the following conditions:
 \begin{itemize}
     \item If $s \leq t$ in $P$, then $\sigma(s) \leq \sigma(t)$; in other words, $\sigma$
is order-preserving.\footnote{Stanley defines $(P, \omega)$-partitions in an order-\emph{reversing} fashion; the present definition mirrors that of the usual order polynomial.}
\item If $s < t$ and $\omega(s) > \omega(t)$, then $\sigma(s) < \sigma(t)$.
 \end{itemize}
Let $\Omega_{P,\omega}(j)$ be the number of $(P, \omega)$-partitions $\sigma: P \rightarrow
[j]_0 := \{ 0, 1, \dots, j \}$.
\end{definition}


The function $\Omega_{P,\omega}(j)$ turns out to be a polynomial of degree $m$, called the
\emph{order polynomial} of $(P,\omega)$.\footnote{
This definition is found in the literature most frequently for the case where $\omega$ is a natural
(i.e., order-preserving) labeling.
} Subsequently, we may define the {\em $h^*$-polynomial} of a labeled poset $(P,\omega)$, $h^*_{ P,\omega }(x)$, via
\[
\sum_{j\geq 0} \Omega_{P,\omega}(j) \, x^j \ = \ \frac{ h^*_{ P,\omega } (x) }{ (1-x)^{m+1} }  \, .
\]
When $\omega$ is a natural labeling, we denote $h^*_{ P,\omega } (x)$ by $h^*_{ P } (x)$.
For example, we can rephrase Corollary~\ref{cor:h*asnarayana} as $h^*_{ [2] \times [n] } (x) = N_n(x)$.

Parallel to the classical case, we may think of a linear extension $\sigma$ of $P$ as a
permutation of $\omega$; we denote the set of all such linear extensions by $\cL(P,\omega)$. 
The fundamental property of order polynomials is the following~\cite[Proposition~13.3]{stanleythesis} (see also~\cite[Theorem~3.15.8]{ec1}).
\begin{theorem}[Stanley]\label{thm: 3.15.8}
$\displaystyle h^*_{P,\omega}(x) = \sum\limits_{\sigma \in \cL(P,\omega)}x^{\des(\sigma)} . $
\end{theorem}

Given a poset $P$ with a labeling $\omega$ and a chain $C$ of $P$, a \emph{descent} of $C$ is any occurrence of a cover relation $a < b$ in $C$ with $\omega(a) > \omega(b)$.
We define $\Des(C,\omega)$ to be the set of descents in $C$ and let $\des(C,\omega)$ denote the cardinality of $\Des(C,\omega)$. The following theorem, which was stated in \cite[Theorem~4.1]{stembridge} using the notion of shift equivalence,  provides a sufficient condition for the $h^*$-polynomials corresponding to two different labelings of a given poset to be the same up to a shift.

\begin{theorem}[Stembridge \cite{stembridge}]\label{thm: main extended}
Let $P$ be a poset with two labelings $\omega$ and $\omega'$ such that for each $j \in P$ there exists $t_j$
with the following conditions:
\begin{itemize}
    \item if $j$ is minimal then $t_j = 0$;
    \item if $j$ covers $i$ then for any maximal chain $C$ containing $i$ and $j$
\[
  t_j - t_i = \begin{cases}
     1 & \text{ if } i\in \Des(C ,\omega) \setminus \Des(C ,\omega') \, , \\
    -1 & \text{ if } i\in \Des(C ,\omega') \setminus \Des(C ,\omega) \, , \\
     0 & \text{ otherwise; }
  \end{cases}
\]
    \item there exists $k$ such that $t_j = k$ for any maximal~$j$.
\end{itemize}
Then $h^*_{P,\omega}(x) = x^k \, h^*_{P,\omega'}(x) \, .$
\end{theorem}

We note that our conditions in Theorem~\ref{thm: main extended} imply that every maximal chain $C$ in the order ideal $\langle j \rangle$
satisfies $t_j=\des(C,\omega)-\des(C,\omega')$.
By choosing $\omega'$ to be natural (and
so there are no descents) and $t_j=\des(C,\omega)$ for any choice
of maximal chain $C$ in $\langle j\rangle $ for a given $j\in P$, we obtain the
following result which appeared in \cite[Corollary~2.4]{branden}, 
phrased there in the language of sign-graded posets.\footnote{The notion that all maximal
chains of $(P,\omega)$ have the same number of descents is slightly different than $(P,\omega)$ being
sign-graded; when $P$ is graded, the two notions coincide.}

\begin{corollary}\label{cor: main}
Consider a poset $P$ with a labeling $\omega$ such that all maximal chains have the same number $k$ of descents. Then $h^*_{P,\omega}(x) = x^k \, h^*_P(x)$. 
\end{corollary}

Our results in the remainder of this section concern posets with maximal chains with the same number
of descents for the sake of simplicity, but most arguments extend to the situation described in Theorem~\ref{thm: main extended}.

\begin{theorem}\label{thm: multiplication}
Consider a poset $P$ with a labeling $\omega$ such that all maximal chains have the same number $k$ of descents. Then $C_n^{P,\omega}(x) = x^k \, A_n(x) \, h^*_{P \times [n]}(x)$.
Moreover, if $P$ is graded then $C_n^{P,\omega}(x)$ is palindromic.
\end{theorem}

Theorem~\ref{thm:main} follows as an immediate corollary, since we have no descents in
$P$ with a natural labeling.

\begin{proof}
For given $\sigma \in S_n$, all maximal chains in $(P \times [n] , \omega \times \sigma)$ have the
same number of descents, namely $k+\des(\sigma)$. Using Corollary~\ref{cor: main}, we compute
    \begin{align*}
        C^{P,\omega}_n(x) &= \sum\limits_{\pi \in \mathcal{C}_n^{P,\omega}}x^{\des(\pi)}=\sum\limits_{\sigma\in
S_n}\sum\limits_{\pi\in \cL(P \times [n] , \omega \times \sigma)} x^{\des(\pi)} 
        \\ &=\sum\limits_{\sigma\in S_n}x^{\des(\sigma)+k} \, h^*_{P \times [n]}(x)=x^k
A_n(x) \, h^*_{P \times [n]}(x) \, . 
        \end{align*}
The $h^*$-polynomial of a naturally labeled poset is palindromic if and only if the poset is
graded~\cite[Proposition~19.3]{stanleythesis} (see also \cite[Corollary~3.15.18]{ec1}).
It follows, in this case, that $C_n^{P,\omega}(x)$ is palindromic.
\end{proof} 

Starting with $P\times [n]$, we now construct a new poset $P \widecheck\times [n]$ 
of cardinality $(m+1)n$ by
adding $n$ elements with no relation among them but  which cover all (maximal)
elements of $P\times [n]$.\footnote{$P \widecheck\times [n]$ is the so-called \emph{ordinal sum} between $P\times [n]$ and an $n$-element antichain.}
 For examples, see the posets in Figure~\ref{fig: labeled poset extended}. 
We extend a given labeling $\omega$ of $P$, first to a labeling $\omega \times \id$ of $P \times
[n]$, and then, to a labeling $\omega \widecheck\times \id$ of $P \widecheck\times [n]$, by giving the new elements any labels that are larger
than those in $\omega \times \id$.

 \begin{figure}[ht]
    \centering
\begin{tikzpicture}[main/.style = {}] 
{\color{magenta}
\node[main] (1) {$1$}; 
\node[main] (2) [above of=1] {$2$};
\node[main] (3) [right of=2] {$3$};
\node[main] (4) [above of=3] {$4$};
\node[main] (5) [right of=4] {$5$};
\node[main] (6) [above of=5] {$6$};
\node[main] (8) [above left of=6] {$7$};
\node[main] (10) [above of=6] {$8$};
\node[main] (11) [above right of=6] {$9$};
\draw (2) -- (1);
\draw (3) -- (4);
\draw (5) -- (6);
\draw (3) -- (1);
\draw (2) -- (4);
\draw (3) -- (5);
\draw (6) -- (4);
\draw (6) -- (8);
\draw (6) -- (10);
\draw (6) -- (11);}
\end{tikzpicture} 
\hspace{1 cm}
\begin{tikzpicture}[main/.style = {}] 
{\color{blue}
\node[main] (1) {$3$}; 
\node[main] (2) [above left of=1] {$2$};
\node[main] (3) [above right of=1] {$1$};
\node[main] (4) [above of=3] {$4$};
\node[main] (3r) [above right of=3] { }; 
\node[main] (5) [below right of=3r] {$7$}; 
\node[main] (6) [above left of=5] {$6$};
\node[main] (7) [above right of=5] {$5$};
\node[main] (8) [above of=7] {$8$};
\node[main] (7r) [above right of=7] { }; 
\node[main] (9) [below right of=7r] {$11$}; 
\node[main] (10) [above left of=9] {$10$};
\node[main] (11) [above right of=9] {$9$};
\node[main] (12) [above of=11] {$12$};

\node[main] (12ul) [above left of=12] { }; 
\node[main] (13) [above left of=12ul] {$13$};
\node[main] (14) [right of=13] {$14$};
\node[main] (15) [ right of=14] {$15$};
\draw (2) -- (1);
\draw (3) -- (1);
\draw (3) -- (4);
\draw (5) -- (6);
\draw (5) -- (7);
\draw (7) -- (8);
\draw (10) -- (9);
\draw (9) -- (11);
\draw (12) -- (11);

\draw (1) -- (5);
\draw (9) -- (5);
\draw (2) -- (6);
\draw (10) -- (6);
\draw (3) -- (7);
\draw (7) -- (11);
\draw (8) -- (12);
\draw (8) -- (4);
\draw (12) -- (9);
\draw (5) -- (9);
\draw (10) -- (13);
\draw (10) -- (14);
\draw (10) -- (15);
\draw (12) -- (13);
\draw (12) -- (14);
\draw (12) -- (15);
}
\end{tikzpicture} 
    \caption{On the \textcolor{magenta}{left}, the naturally labeled poset for $(P \widecheck\times [3],\omega \widecheck\times \id)$ and $P=[2]\times[3]$ satisfying $h^*_{P
\widecheck\times [3],\omega \widecheck\times \id}(x)=C_3^2(x)$.
\\  { On the {\textcolor{blue}{right}, the labeled poset $(P \widecheck\times [3],\omega \widecheck\times \id)$ for  $P=T\times [3]$ and $(T,\omega)$}  labeled poset on $4$ nodes 
satisfying  $h^*_{P \widecheck\times [3],\omega \widecheck\times \id}(x)=C_3^{P,\omega}(x)$.}
}
    \label{fig: labeled poset extended}
\end{figure}

\begin{proposition}\label{prop:oneposet} 
\ $ C_n^{P,\omega}(x)=h^*_{P \widecheck\times [n],\omega \widecheck\times \id}(x) \, . $
\end{proposition}

We remark that this gives a direct proof for the palindromicity of $C_n^{P,\omega}$ when $\omega=\id$ since $P
\widecheck\times [n]$ is graded. 

\begin{proof} 
By construction of $(P \widecheck\times [n],\omega \widecheck\times \id)$, any linear extension of $(P \widecheck\times [n],\omega \widecheck\times \id)$ is the concatenation of a
linear extension of $(P\times [n],\omega \times \id)$ and a linear extension of an $n$-element
antichain---that is, a permutation of $S_n$.
Moreover, since by construction there is no descent in position $nm$ in any linear
extension $\pi\in \mathcal{L}(P \widecheck\times [n],\omega \widecheck\times \id)$, the number of descents of $\pi$ is the sum of the
number of descents in the two linear extensions concatenated. 
This describes a bijection between linear extensions of $(P \widecheck\times [n],\omega \widecheck\times \id))$ and pairs of a linear
extension of $(P\times [n], \omega \times \id)$ and a permutation of $S_n$. 
It follows that \[\sum_{\pi\in \mathcal{L}(P \widecheck\times [n],\omega \widecheck\times \id) }
x^{\des(\pi)}=\sum_{\pi\in \mathcal{L}(P\times [n],\omega \times \id)}
x^{\des(\pi)}\sum_{\sigma\in S_n}
x^{\des(\sigma)}.\]

On the other hand, by applying Corollary~\ref{cor: main} to the poset $P
\times [n]$ and labelings $\omega\times \sigma$ and $\omega\times \id$,
we see that the corresponding $h^*$-polynomials are the same up to a shift, i.e., 
\[\sum_{\pi\in \mathcal{L}(P\times [n],\omega \times \sigma)}
x^{\des(\pi)}=x^{\des(\sigma)}\sum_{\pi\in \mathcal{L}(P\times [n],\omega \times \id)} 
x^{\des(\pi)}.\] Hence
\[C_n^{P,\omega}(x)=\sum_{\sigma \in S_n}\sum_{\pi\in \mathcal{L}(P\times [n],\omega \times \sigma)}
x^{\des(\pi)}=\sum_{\sigma \in S_n}x^{\des(\sigma)}\sum_{\pi\in \mathcal{L}(P\times [n],\omega
\times \id)} x^{\des(\pi)}=h^*_{P \widecheck\times [n],\omega \widecheck\times \id}(x) \, . \qedhere \] 
\end{proof}

\begin{remk}
    One may notice that the product factorization structure holds more generally (and so do the palindromicity results and the existence of a poset similar to the one in Proposition~\ref{prop:oneposet}). In particular, we may consider $(P,\omega)$ with exactly $k$ descents across each of its maximal chains. Let 
  $P'$ be a poset on $n$ elements. Then it is a consequence of Theorem~\ref{thm: main extended} that

  \begin{align*}
        \sum\limits_{\sigma \in \mathcal{L}(P')}
        h^*_{P \times [n] , \omega \times \sigma} (x) &=x^{k} \, h^*_{P'}(x) \, h^*_{P \times [n]}(x) \, .
        \end{align*}
        When $P'$ is the antichain, $h^*_{P'}=A_n$ and we recover the identity in Theorem~\ref{thm: multiplication}.
\end{remk}


\section{Amphibians}\label{sec:amphibians}

In this section, we will study a broader family of linear extensions/multiset permutations
whose descent polynomials are also palindromic.

Throughout this section,  we denote by $Q$ a subposet of $P\times [n]$ with some relations missing between different
copies of $P$, i.e., those of the form $(p,j) \prec (p,j+1)$. 
An example of such a poset $Q$, with two different labelings, is shown in Figure~\ref{fig:labeledsubposet}.
The motivation for studying these subposets comes from the following observation. When $P=[m]$, the linear
extensions of $Q$ can be interpreted (in the same way as before) as a collection of multiset
permutations of $\{1^m,2^m,\dots,n^m\}$, where conditions weaker than those for canon
permutations are imposed. 
When $P=[m]$ and we remove all relations of the form $(p,j) \prec (p,j+1)$ for all  $p\in [m-1]$, 
 we recover the set of multiset permutations of
$\{1^m,2^m,\dots,n^m\}$ whose subsequence formed by the  last copy of each element in $[n]$ is fixed. Similarly to canon permutations, summing over the corresponding descent polynomials for all canon labelings $\id \times \sigma$  of $Q$ for $\sigma \in S_n$ gives the descent polynomial of multiset permutations of $\{1^m,2^m,\dots,n^m\}$. 
Since the descent polynomials corresponding to $\{1^m,2^m,\dots,n^m\}$ and $\mathcal{C}_n^m$
are both palindromic, it is natural to ask if palindromicity extends to more subposets $Q$ of
 $P\times [n]$.
Theorem~\ref{them: subposet} below confirms this. 
We note that the notion of \emph{palindromicity} here might involve polynomials with zero constant
terms (or even more zero coefficients), and so we state, in each case, the relevant functional
equation.

Let $Q$ be a subposet of $P\times [n]$ with some of the cover relations of the form $(p,j) \prec (p,j+1)$ removed.
We define the \emph{dissonant canon polynomial} 
\[ C^{Q,\omega}(x):=\sum_{\sigma \in S_n} h^*_{Q,{\omega}\times \sigma}(x) \, .\]
The name is inspired from the case when $P=[m]$ and $Q$ is a subposet of $[m]\times [n]$ with some of
the edges of the form $(p,j) \prec (p,j+1)$ removed, where $p\in P\setminus\{q\}$ for some fixed
$q\in P$. 
 
This condition fixes one of the occurrence subsequences,
so the union over \(\sigma\in S_n\) is disjoint and \(C^{Q,\id}(x)\) is the
ordinary descent polynomial of the corresponding multiset permutations.
 We call this family of multiset permutations of $\{1^m, \dots,n^m\}$ the \emph{dissonant canon permutations} (for the corresponding $Q$).



If all maximal chains in $(P,\omega)$ have the same number
$k$ of descents, then any maximal chain in $(Q, \omega\times \sigma)$ will contain between $k$ and
$k+\des(\sigma)$ descents. In particular, from Theorem~\ref{thm: main extended} we can deduce the following.
\begin{corollary}
    \label{cor: main extended}
Consider a poset $P$ with a labeling $\omega$ such that all maximal chains in $(P,\omega)$ have the same number $k$ of descents. 
Then, for any $Q$ as above 
\[h^*_{Q,\omega\times \sigma}(x) = x^k \, h^*_{Q,\id\times \sigma}(x) \, .\]
\end{corollary}

\begin{proof}
For $p\in P$, let $s_p$ be the number of descents in any maximal chain of
the principal order ideal $\langle p\rangle$ with respect to $\omega$, and
set
\[
t_{(p,i)}:=s_p 
\]
for $(p,i)\in Q$. 
This is well-defined by the assumption on maximal chains of $(P,\omega)$,
since two maximal chains in $\langle p\rangle$ can be extended by the
same chain above $p$. 
Our goal is to verify the hypotheses of Theorem~\ref{thm: main extended} for the two labelings
$\omega\times\sigma$ and $\operatorname{id}\times\sigma$ of $Q$, using the
above choice of $t_{(p,i)}$. 

For a cover of the form $(p,i)\prec (p',i)$, we have
\[
t_{(p',i)}-t_{(p,i)}=s_{p'}-s_p \, .
\]
This difference is equal to $1$ precisely when $p\prec p'$ is a descent
with respect to $\omega$, and is equal to $0$ otherwise, because a maximal
chain in $\langle p\rangle$ together with $p'$ is a maximal chain in
$\langle p'\rangle$. Since $\operatorname{id}$ has no descents on $P$, this
is exactly the required local difference between the labelings
$\omega\times\sigma$ and $\operatorname{id}\times\sigma$. 
On the other hand, for a cover of the form $(p,i)\prec (p,i+1)$, the two labelings have a
descent simultaneously, namely precisely when $\sigma(i)>\sigma(i+1)$.
Thus the local difference is $0$, which agrees with
\[
t_{(p,i+1)}-t_{(p,i)}=0 \, .
\]

Lastly, if $(p,i)$ is minimal in $Q$, then $p$ is minimal in $P$, so $s_p=0$. If $(p,i)$ is maximal in $Q$, then $p$ is maximal in $P$, so a maximal chain in $\langle p\rangle$ is a maximal chain of $P$ and hence has $k$ descents. 
Therefore the hypotheses of Theorem~\ref{thm: main extended} are satisfied.
\end{proof}

Corollary \ref{cor: main extended} will be used to prove the palindromicity of the polynomials considered in this section.


\begin{theorem}\label{them: subposet}
Consider a graded,  $m$-element poset $P$ with a labeling $\omega$ such that all maximal chains in $(P,\omega)$ have the same number $k$ of descents (and $k'$ ascents).
Then, for any $Q$ as above, the dissonant canon polynomial $C^{Q,\omega}(x)$ is palindromic, 
in the sense that
\[
 x^{mn-1+k-k'} \, C^{Q,\omega}(\tfrac{1}{x})= C^{Q,\omega}({x}) \, .
\]
\end{theorem}

\begin{proof}
We define the bijection $\phi: S_n\rightarrow S_n$ as follows. 
Given $\sigma \in S_n$, let
$\phi(\sigma)$ be the permutation of  $S_n$ where the entry $j$ is replaced by $n+1-j$, for all $j\in [n]$.  For example, if $\sigma=1342$ then $\phi(\sigma) = 4213$.
We extend $\phi$ to a map $\phi: \mathcal{L}(Q,\omega\times\sigma)\rightarrow \mathcal{L}(Q,\phi({\omega})\times {\phi(\sigma)})$, sending $(\omega\times \sigma)(p,j)$ to $(\phi({\omega})\times \phi(\sigma))(p,j)$ for $(p,j)\in Q$, where $\phi({\omega}(p))=m+1-\omega(p)$.

For every \((p,j)\in Q\), 
\[
(\omega\times\sigma)(p,j)
+
(\phi(\omega)\times\phi(\sigma))(p,j)
=mn+1 \, .
\]
Thus \(\phi(\omega)\times\phi(\sigma)\) is obtained from
\(\omega\times\sigma\) by reversing all labels. 
Hence the map $\phi$ converts descents of
$\pi\in\mathcal{L}(Q,\omega\times\sigma)$ to ascents of
$\phi(\pi)\in\mathcal{L}(Q,\phi(\omega)\times\phi(\sigma))$ and vice versa,
since $\phi$ is an involution.

Therefore, since $\text{asc}(\pi)=mn-1-\des(\pi)$ for $\pi\in \mathcal{L}(Q,\phi({\omega})\times {\phi(\sigma)})$,

\begin{equation}\label{eq:reciprocityidentity}h^*_{Q,\omega\times\sigma} \left(\frac{1}{x} \right)=
\sum\limits _{\pi\in \mathcal{L}(Q,\phi({\omega})\times {\phi(\sigma)})} \left(\frac{1}{x} \right)^{\text{asc}(\pi)}
=
\frac{1}{x^{mn-1}}
\sum\limits _{\pi\in \mathcal{L}(Q,\phi({\omega})\times {\phi(\sigma)})}{x^{\des(\pi)}}.
    \end{equation}
By construction of $Q$, every maximal chain in $(Q,\phi({\omega})\times {\phi(\sigma)})$ (resp.
$(Q,{\omega}\times {\phi(\sigma)})$) has $\des(\phi(\omega))$ (resp. $\des(\omega)$) more descents than the corresponding maximal chain in $(Q,\id\times {\phi(\sigma)})$.\footnote{ 
Here we use the notation $\des(\omega)=k$ if each maximal chain in $(P,\omega)$ has $k$ descents.
The fact that $\des(\phi(\omega))$ is well-defined follows by the graded assumption.}
We now apply Corollary~\ref{cor: main extended}, yielding
\begin{equation}\label{eq:applyingmainlemma}
\sum\limits _{\pi\in \mathcal{L}(Q,\phi({\omega})\times {\phi(\sigma)})}{x^{\des(\pi)}}=x^{\des(\phi({\omega}))-\des(\omega)}
\sum\limits _{\pi\in \mathcal{L}(Q,{\omega}\times {\phi(\sigma)})}{x^{\des(\pi)}}.
\end{equation}

Equations \eqref{eq:reciprocityidentity} and \eqref{eq:applyingmainlemma} combine into
\begin{equation}\label{eq: bijection refinement}
    h^*_{Q,\omega\times\sigma} \left(\frac{1}{x} \right)=\frac{x^{\des(\phi(\omega))-\des(\omega)}}{x^{mn-1}}
\sum\limits _{\pi\in \mathcal{L}(Q,{\omega}\times {\phi(\sigma)})}{x^{\des(\pi)}}
=\frac{x^{\des(\phi(\omega))-\des(\omega)}}{x^{mn-1}} \, h^*_{Q,{\omega}\times {\phi(\sigma)}}{(x)} \, .
\end{equation}
Since $\phi:S_n \rightarrow S_n$ is a bijection, it follows that 
\begin{equation*}
    \sum_{\sigma\in S_n}h^*_{Q,\omega\times\sigma} \left(\frac{1}{x} \right)=
\frac{x^{\des(\phi({\omega}))-\des(\omega)}}{x^{mn-1}}
\sum_{\sigma\in S_n}h^*_{Q,\omega\times\sigma}({x}) \, . \qedhere
\end{equation*}

\end{proof}

A key point in the proof of Theorem~\ref{them: subposet} is the decomposition of the polynomial whose
palindromicity is in question into (a total of $\frac{n!}{2}$) palindromic polynomials (see Equation \eqref{eq: bijection refinement}). A natural question is to identify the corresponding bijections between the permutations captured by the coefficients, which would result in a bijective proof of Theorem~\ref{them: subposet}. This would also yield bijections in the subcase of canon permutations, addressing questions raised by Elizalde \cite[Problem~4.1]{sergi}.

\begin{corollary}\label{cor: palindromicity direct cons}
    Consider a chain poset $P=[m]$ with a labeling $\omega$ whose linear extension has $k$ descents. 
Then, for any $Q$ as above, the dissonant canon polynomial $C^{Q,\omega}(x)$ satisfies 
\[
 x^{m(n-1)+2k}  \, C^{Q,\omega}(\tfrac{1}{x})= C^{Q,\omega}({x}) \, ,
\] where $m(n-1)+k$ is the degree of $C^{Q,\omega}(x)$.
\end{corollary}

\begin{proof}
The palindromicity follows from Theorem \ref{them: subposet} by substituting $k'=\des(\phi({\omega}))=m-1-\des(\omega)=m-1-k$.

Let us now compute the degree of $C^{Q,\omega}(x)$. 
Consider the elements in $P\times[n]$ of the form
$(p,\sigma^{-1}(1))$ with $p\in P$. By construction of the labeling
$\omega\times\sigma$, these are precisely the elements whose labels lie in
$[m]$. Hence any such element can form a descent in a linear extension of
$(Q,\omega\times\sigma)$ only if it is followed by another element
$(p',\sigma^{-1}(1))$ such that $\omega(p)>\omega(p')$.
Therefore, at least $m-k$ of the elements $(p,\sigma^{-1}(1))$ in $Q$ cannot occur as the first entry of a descent in any linear extension of $(Q,\omega\times\sigma)$. 
Hence the degree of $h^*_{Q,\omega\times \sigma} $ is bounded above by $m(n-1)+k$, for any $\sigma
\in S_n$. 

For the lower bound, let
\(\pi=\pi_1\pi_2\cdots\pi_m\) be a linear extension of \((P,\omega)\)
with \(k\) descents. Write \(\pi\) as a concatenation of its maximal
descending runs,
\[
\pi=\tau^{(1)}\tau^{(2)}\cdots\tau^{(r)},
\]
where if \(\tau^{(s)}=\pi_a\pi_{a+1}\cdots\pi_b\), then
\(\pi_i>\pi_{i+1}\) for \(a\leq i<b\), and the boundary between two
consecutive runs is not a descent. Then
\[
\sum_{s=1}^r (|\tau^{(s)}|-1)=k
\qquad\text{and}\qquad
r=m-k \, .
\]
For \(\tau^{(s)}=\pi_a\pi_{a+1}\cdots\pi_b\), we let 
\[
\tau^{(s)}+tm=(\pi_a+tm)(\pi_{a+1}+tm)\cdots(\pi_b+tm) \, .
\]
Consider the word
\[
\pi_{\text{max}}=
\prod_{s=1}^r
\bigl(\tau^{(s)}+(n-1)m\bigr)
\bigl(\tau^{(s)}+(n-2)m\bigr)
\cdots
\bigl(\tau^{(s)}+m\bigr)
\tau^{(s)} ,
\]
where the product denotes concatenation in the order \(s=1,\ldots,r\). The permutation \(\pi_{\text{max}}\) is a linear extension of \((Q,\omega\times w_0)\), where \(w_0\) is the decreasing permutation, \(w_0(j)=n+1-j\).
If
\(\tau^{(s)}\) has length \(\ell_s\), then its \(n\) copies contribute
\(n\ell_s-1\) descents, where the count includes 
 $n-1$ boundary descents (coming from the concatenation of the copies of \(\tau^{(s)}\)). Thus, the total 
number of descents of \(\pi_{\text{max}}\) is
\[
\sum_{s=1}^r (n\ell_s-1)
=
nm-r
=
m(n-1)+k \, .
\]
Hence \(\pi_{\max}\) attains the upper bound, which completes the proof.
\end{proof}

The palindromicity of the dissonant canon polynomial
$C^{Q, \id}(x)$ in Theorem~\ref{them: subposet} also gives rise to the following problem. 

\begin{question}\label{que: gampos}
    For which subposets $Q$ of $P{\times}[n]$ is $C^{Q,\id}(x)$ $\gamma$-positive? Can we describe the $\gamma$-coefficients in a unified way?
\end{question}

In the next section (Corollary~\ref{cor:gamma-interp-canon}), we will see that the answer to Question~\ref{que: gampos} is positive for canon permutations, which correspond to the poset $Q=[m]\times [n]$.
The answer is also positive for the poset corresponding to all multiset
permutations, i.e., the poset where all conditions of the form $(p,j)\prec (p,j+1)$ are removed. 
This follows from \cite{branden} by realizing the set of multiset permutations of $\{1^m,2^m,\dots,n^m\}$ as linear extensions of the poset consisting of $n$ $m$-element chains, labeled with the regular canon labeling.
A combinatorial interpretation for those $\gamma$-coefficients is given in \cite{BrandenActions}.

The two combinatorial interpretations (corresponding to the set of multiset permutations and to canon
permutations) rely on partitioning permutations into classes, but the classes are different
in each case. Question~\ref{que: gampos} asks for a unified interpretation of the
$\gamma$-coefficients for both corresponding polynomials, as well as other $\gamma$-positive
polynomials that may arise for different choices of subposets $Q$.  We conjecture that this class will at least contain the descent polynomials of dissonant canon permutations.

 Let $Q$ be a subposet of $[m]\times [n]$ with some of the cover relations of the form $(p,j) \prec (p,j+1)$ removed, for $p\in [m-1]$.
Then \(C^{Q,\id}(x)\) is the descent polynomial of the corresponding class of
dissonant canon permutations of \(\{1^m,2^m,\ldots,n^m\}\).

\begin{conj}\label{conj: gammapos}
 
The descent polynomial $C^{Q,\id}(x)$ is $\gamma$-positive, for any class of dissonant canon permutations.
\end{conj}

We have confirmed Conjecture~\ref{conj: gammapos}  for the pairs \[(m,n)\in \big\{(2,2),\,(2,3),\,(2,4),\,(2,5),\,(2,6),\,(3,2),\,(3,3),\,(3,4),\,(3,5),\,(4,2),\,(4,3),\,(4,4)\big\}.\]
However, the method used in the next section for the canon permutation case could not be applied, since we are
not aware of a single poset encoding the collection of all labeled linear extensions in
$C^{Q,\id}(x)$ for general $Q$. 

\begin{remk}
    In the proof of Theorem~\ref{them: subposet}, we decomposed $C^{Q,\id}$ into a sum of palindromic
polynomials. We note that these summands are not always $\gamma$-positive, even for the canon posets.
For example, for $\phi(\id)= 321$,
\[
h^*_{[2]\times [3],\id \times \id}(x) +h^*_{[2]\times [3],\id \times \phi(\id)}(x)=(1+3x+x^2)(1+x^2)\, ,
\] 
which is palindromic but not $\gamma$-positive---in fact, its coefficient vector is not even unimodal. 
However, we believe that $\gamma$-positivity holds for $C^{Q,\id}(x)$ (as it does for canon polynomials). 
\end{remk}

\begin{figure}[ht]
    \centering
\begin{tikzpicture}[main/.style = {}] 
{\color{magenta}
\node[main] (1) {$1$}; 
\node[main] (2) [above of=1] {$2$};
\node[main] (3) [right of=2] {$5$};
\node[main] (4) [above of=3] {$6$};
\node[main] (5) [right of=4] {$7$};
\node[main] (6) [above of=5] {$8$};
\node[main] (7) [right of=6] {$3$};
\node[main] (8) [above of=7] {$4$};
\draw (2) -- (1);
\draw (3) -- (4);
\draw (5) -- (6);
\draw (7) -- (8);
\draw (3) -- (1);
\draw (2) -- (4);
\draw (5) -- (7);
\draw (3) -- (5);
\draw (6) -- (4);}
\end{tikzpicture} 
\begin{tikzpicture}[main/.style = {}] 
{\color{blue}
\node[main] (4) {$8$}; 
\node[main] (2) [above of=1] {$7$};
\node[main] (3) [right of=2] {$4$};
\node[main] (4) [above of=3] {$3$};
\node[main] (5) [right of=4] {$2$};
\node[main] (6) [above of=5] {$1$};
\node[main] (7) [right of=6] {$6$};
\node[main] (8) [above of=7] {$5$};
\draw (2) -- (1);
\draw (3) -- (4);
\draw (5) -- (6);
\draw (7) -- (8);
\draw (3) -- (1);
\draw (2) -- (4);
\draw (5) -- (7);
\draw (3) -- (5);
\draw (6) -- (4);}
\end{tikzpicture} 
    \caption{A subposet $Q$ of $([2]\times [4],\id\times {1342})$ \textcolor{magenta}{(left)} and its corresponding labeling through $\phi$ \textcolor{blue}{(right)}.}
    \label{fig:labeledsubposet}
\end{figure}


\section{$\gamma$-positivity of Canon Permutations}\label{sec: canon positivity}



Theorem~\ref{thm:main} shows that certain distributional properties shared among the polynomials
$A_n(x)$ and $h^*_{P\times [n]}(x)$, such as palindromicity and $\gamma$-positivity, transfer to
the canon polynomial $C_n^{P,\omega}(x)$. 
A direct proof of the $\gamma$-positivity of $C_n^m(x)$ can be derived from Br\"and\'en's work in \cite{branden} using the poset $[m] \widecheck\times [n]$ constructed in Proposition~\ref{prop:oneposet} (here $P=[m]$). 
Below, we discuss a combinatorial interpretation of the $\gamma$-coefficients 
of $C_n^m(x)$  in terms of canon permutations, a consequence of a group action on permutations due to Foata and Strehl extended to posets by Br\"and\'en~\cite{BrandenActions}.

Let $N=(m+1)n$.
For a permutation $\pi=\pi_1\cdots \pi_{N}\in \mathcal{L}([m] \widecheck\times [n])$, we construct the signed permutation 
$\pi_1^{c_1}\cdots \pi_{N}^{c_{N}}$ where $c_i=0$ if the maximal chains in the poset ideal $\langle \pi_i \rangle$ have even length and $c_i=1$ if the length is odd. We also set $\pi_0^{c_0}=\pi_{N+1}^{c_{N+1}}=0^1$. 
We define a \emph{signed descent} of $\pi$ to be a descent of the signed permutation $0^1\pi_1^{c_1}\cdots \pi_{N}^{c_{N}}0^1$, i.e., a position $i\in \{0,1,\ldots, N\}$ such that $\pi_i^{c_i}>\pi_{i+1}^{c_{i+1}}$ according to the ordering $1^0<2^0<\cdots <N^0<0^1<1^1<2^1<\cdots<N^1$. We denote the number of nonboundary signed descents of $\pi\in \mathcal{L}([m] \widecheck\times [n])$ by \(\des_{\pm}(\pi):=\left |\{i\in [N-1]: \pi_i^{c_i}>\pi_{i+1}^{c_{i+1}}\}\right|\). 
We say that \(\pi\) has a \emph{double signed descent} if there exists
\(i\in\{0,1,\ldots,N-1\}\) such that both \(i\) and \(i+1\) are signed descents
of the augmented permutation
$0^1\pi_1^{c_1}\cdots \pi_N^{c_N}0^1$.

\begin{corollary}\label{cor:gamma-interp-canon}
Writing
\[
C_n^m(x)
=
\sum_{i=0}^{\lfloor m(n-1)/2\rfloor}
\gamma_i x^i(1+x)^{m(n-1)-2i},
\]
the coefficient \(\gamma_i\) counts linear extensions
\(\pi \in\mathcal L([m]\widecheck\times[n])\),  such that
\begin{itemize}
\item \(\des_{\pm}(\pi)=i+\lfloor (m+n-1)/2\rfloor\);
\item \(\pi\) has no double signed descents.
\end{itemize}
\end{corollary}
\begin{proof}
By Proposition~\ref{prop:oneposet}, we have
\(C_n^m(x)=h^*_{[m]\widecheck\times[n]}(x)\).
The stated interpretation is then the specialization of Br\"and\'en's formula
for sign-graded posets~\cite[Section~6]{BrandenActions} to the canonical
labeling described above.
\end{proof}

\begin{example}\label{ex: main theorem}


For \(m=3\) and \(n=2\), using the poset $[3] \widecheck\times [2]$ shown in Figure~\ref{fig: example one poset}, we compute 
   \[C_2^3(x)=A_2(x) \, N_3(x) \,=(1+x)(1+3x+x^2) =(1+x)^3+x(1+x) \] with the following combinatorial interpretation by Corollary~\ref{cor:gamma-interp-canon}.

The coefficient of $(1+x)^3$ corresponds to  the linear extension $\pi_1\pi_2\cdots \pi_8=12435678 \in \mathcal{L}(
[3]\widecheck \times [2])$ (that is, the multiset permutation $112122$), whose {signed} counterpart equals
$0^1\,1^0\,2^1\,4^1\,3^0\,5^0\,6^1\,7^0\,8^0\,0^1$. Its nonboundary signed descents are \(3\) and \(6\), whence we
have
\(\des_{\pm}(12435678)=2\), and its 
full set of signed descents is \(\{0,\,3,\,6\}\). It has no double signed descents, and hence contributes to~\(\gamma_0\).

The coefficient of $x(1+x)$ represents the linear extension $\pi_1\pi_2\cdots \pi_8=12345678\in \mathcal{L}([3]\widecheck \times [2])$ (i.e., the multiset permutation $111222$), with corresponding signed permutation \[
0^1\,1^0\,2^1\,3^0\,4^1\,5^0\,6^1\,7^0\,8^0\,0^1.
\]
Its nonboundary signed descents are \(2\), \(4\) and \(6\), so
\(\des_{\pm}(12345678)=3\), and its full set of signed descents is \(\{0,\,2,\,4,\,6\}\). Since there are no double signed descents, it contributes to the coefficient \(\gamma_1\).
\end{example}
\begin{figure}
\begin{tikzpicture}[main/.style = {}] {\color{blue}
\node[main] (1) {$1$}; 
\node[main] (2) [above of=1] {$2$};
\node[main] (3) [above of=2] {$3$};
\node[main] (4) [above right of=1] {$4$};
\node[main] (5) [above of=4] {$5$};
\node[main] (6) [above of=5] {$6$};
\node[main] (7) [above left of=6] {$7$};
\node[main] (8) [above right of=6] {$8$};
\draw (2) -- (1);
\draw (3) -- (2);
\draw (5) -- (6);
\draw (5) -- (4);
\draw (1) -- (4);
\draw (2) -- (5);
\draw (6) -- (3);
\draw (6) -- (8);
\draw (6) -- (7);}
\end{tikzpicture} 
    \caption{The poset $[3] \widecheck\times [2]$.}
    \label{fig: example one poset}
\end{figure}

 As a second example, we describe the $\gamma$-interpretation of $C_3^3(x)$.

\begin{example}
  The canon polynomial corresponding to $m=n=3$ has $\gamma$-vector $(1,8,14,4)$. Specifically, the canon
permutation contributing to the coefficient $\gamma_0=1$ is
    $112123233$ (after translating via linear extensions).
    The permutations contributing to $\gamma_1=8$ are
    \(111223233,\) \(112132233,\) \(112231233,\) \(112122333,\) \(112233123,\)
    \(123112233,\)
    \(221213133,\) \(331312122.\)
   The permutations contributing to $\gamma_2=14$ are
   \(111222333\), \(123123123\),
   \(222113133\), \(221231133\), \(221132133\), \(221211333\), \(221133213\), \(213221133\), 
   \(333112122\), \(331321122\), \(331123122\), \(331311222\), 
   \(331122312\), \(312331122\).
   Lastly, for $\gamma_3=4$ we have the permutations
    \(222111333,\)  \(333111222,\)  \(213213213,\)  \(312312312\).
Each of these permutations represents a class described in \cite{BrandenActions}, and these classes partition the space of canon permutations.
\end{example}

Let us now discuss the polynomial corresponding to weak descents of canon permutations.

\begin{proposition}\label{prop:weakdes} 
$ \displaystyle \sum\limits_{\pi \in\mathcal{C}_n^m}x^{\operatorname{wdes}(\pi)}=
x^{m-1}C^{m}_n(x)=x^{m-1}A_{n}(x) \, h^*_{[m]\times [n]}(x) \, . $
\end{proposition}

\begin{proof}
    As discussed, we can identify the linear extensions of
$([m]\times [n],\omega\times\sigma)$ with the canon permutations in $\mathcal{C}_n^{m,\sigma}$
 when $\omega = \id$ is the identity (i.e., $[m]$ is naturally labeled), for any $\sigma\in S_n$. 
In general, 
the elements of the $j$th copy of $[m]$ in a linear extension of $([m]\times [n],\omega\times \sigma)$
can be identified with the entries with value $j$ in the corresponding canon permutation in
$\mathcal{C}_n^{m,\sigma}$.
When $\omega=\id$, 
the descents in the linear extensions correspond to the descents in the associated canon permutations.

  Let $w_0$ be the reverse natural labeling  of \([m]\), that is,
\(w_0(i)=m+1-i\). 
A linear extension $\pi=\pi_1\pi_{2}\cdots\pi_{mn}\in \mathcal{C}_n^{[m],{w_0}} = \left\{ \cC_n^{ [m], {w_0} \times \sigma } : \, \sigma \in S_n \right\}$ will have a descent $j\in [mn-1]$ whenever either $\pi_j$ comes from the $i$-th copy of $[m]$, $\pi_{j+1}$ comes from the $(i+1)$-th copy of $[m]$, and $\sigma(i)>\sigma(i+1)$, or both $\pi_j$ and $\pi_{j+1}$ come from the same copy of $[m]$.
This means that descents in the linear extensions of $\mathcal{C}_n^{[m],{w_0}}$ correspond to weak descents in the associated canon permutations. Therefore,
\[\sum\limits_{\pi \in\mathcal{C}_n^{m}}x^{\operatorname{wdes}(\pi)}= C^{[m],{w_0}}_n(x) \, .\]
   Since every maximal chain of \(([m],w_0)\) has \(m-1\) descents, Corollary~\ref{cor: main}
gives
\[C^{[m],{w_0}}_n(x)=x^{m-1}C^{[m],\id}_n(x)=x^{m-1}C^{m}_n(x) \, .\, \qedhere \]
\end{proof}

 It follows that the weak-descent polynomial of canon permutations is palindromic and $\gamma$-positive with a specified combinatorial interpretation arising by shifting the coefficients in Corollary~\ref{cor:gamma-interp-canon}. The palindromicity of the weak-descent polynomial of canon permutations for the multiset $\{1,1,2,2,...,n,n\}$ was first observed in   \cite[Corollary~2.3]{sergi}.


  The same argument applies to the dissonant canon permutations from
Section~\ref{sec:amphibians}. Combining Proposition~\ref{prop:weakdes} with
Corollary~\ref{cor: palindromicity direct cons} gives the following.


\begin{corollary}\label{prop:palindromic subposet}
  For a class of dissonant canon permutations of a poset $Q$, the corresponding descent polynomial $C^{Q, \id}(x)$ is palindromic in the sense that
\[
 x^{m(n-1)} \, C^{Q,\id}(\tfrac{1}{x})= C^{Q,\id}({x}) \, .
\]

The dissonant canon polynomial $C^{Q, {w_0}}(x)$ corresponding to weak descents of dissonant canon permutations is palindromic, in the sense that 
\[
 x^{m(n+1)-2} \, C^{Q,w_0}(\tfrac{1}{x})= C^{Q,w_0}({x}) \, .
\] 
\end{corollary} 

  This extends the palindromicity phenomenon from canon permutations to dissonant canon permutations, for both descents and weak descents.

\section*{Acknowledgements}

\thanks{We thank Petter Br\"and\'en and Ben Braun for helpful pointers to the literature.  We are grateful to anonymous referees for their careful review and many helpful suggestions. D.\ D.\ acknowledges support from the Wallenberg Autonomous Systems and
Software Program (WASP), funded by the Knut and Alice Wallenberg Foundation, by the Spanish--German
project COMPOTE (AEI $PCI2024-155081-2$ and DFG $541393733$), and by the
Spanish project PID$2022-137283$NB$\-\ $C$21$ of MCIN/AEI/$10.13039/501100011033$.  A large part of this work took place during a research visit of D.\ D.\ to San Francisco State University in the Fall of 2023, funded by WASP.}

\bibliographystyle{amsplain}
\providecommand{\bysame}{\leavevmode\hbox to3em{\hrulefill}\thinspace}
\providecommand{\MR}{\relax\ifhmode\unskip\space\fi MR }
\providecommand{\MRhref}[2]{%
  \href{http://www.ams.org/mathscinet-getitem?mr=#1}{#2}
}
\providecommand{\href}[2]{#2}

\end{document}